\documentclass{article}
\usepackage{graphicx}
\usepackage{amsmath}
\usepackage{amsfonts}
\usepackage{amsthm}

\theoremstyle{definition}

\begin{document}

\title{Number Gossip}

\author{Tanya Khovanova \\
\textit{Department of Mathematics, MIT}}

\maketitle

\abstract{This article covers my talk at the Gathering for Gardner 2008, with some additions.}


\section{Introduction}

My pet project Number Gossip has its own website: http://www.numbergossip.com/, where you can plug in your favorite integer up to 9,999 and learn its properties. A behind-the-scenes program checks your number for 49 regular properties and also checks a database for unique properties I collected.

\section{Eight}

The favorite number of this year's Gathering is composite, deficient, even, odious, palindromic, powerful, practical and Ulam. It is also very cool as it has the rare properties of being a cake and a narcissistic number, as well as a cube and a Fibonacci number. And it also is a power of two.

In addition, eight has the following unique properties:

\begin{itemize}
\item 8 is the only composite cube in the Fibonacci sequence
\item 8 is the dimension of the octonions and is the highest possible dimension of a normed division algebra
\item 8 is the smallest number (except 1) which is equal to the sum of the digits of its cube
\end{itemize}

\section{Properties}

There are 49 regular properties that I check for:

\begin{tabular}{ l l l l l }
  abundant & evil & odious & Smith \\
  amicable & factorial & palindrome & sociable \\
  apocalyptic power & Fibonacci & palindromic prime & square \\
  aspiring & Google & pentagonal & square-free \\
  automorphic & happy & perfect & tetrahedral \\
  cake & hungry & power of 2 & triangular \\
  Carmichael & lazy caterer & powerful & twin \\
  Catalan & lucky & practical & Ulam \\
  composite & Mersenne & prime & undulating \\
  compositorial & Mersenne prime & primorial & untouchable \\
  cube & narcissistic & pronic & vampire \\
  deficient & odd & repunit & weird \\
  even &  &  &  \\
\end{tabular}

I selected regular properties for their importance as well as their funny names, so your favorite number could be lucky and happy at the same time, as is the case for 7.

Here are some definitions:

The number $n$ is called an \textbf{apocalyptic power} if $2^n$ contains 666. The smallest apocalyptic power is 157: $2^{157} = 182687704666362864775460604089535377456991567872$. The bigger your number the more likely it is to be an apocalyptic power. As your numbers grow the apocalypse is becomes more and more probable.

The number $n$ is \textbf{evil} if it has an even number of ones in its binary expansion. Can you guess what \textbf{odious} numbers are (see John Conway at al. \cite{Conway})?

The $n^{th}$ \textbf{lazy caterer} number is the maximum number of pieces a circular pizza can be cut into with $n$ straight-line cuts. For some reason mathematicians think that pizza is a 2-dimensional object. That means that horizontal cuts are not allowed. A 3D pizza is called a cake, so you can deduce which are \textbf{cake} numbers. For cake numbers, horizontal cuts are allowed, but beware: someone may be deprived of icing.

The $k^{th}$ \textbf{hungry} number is the smallest number $n$ such that $2^n$ contains the first $k$ digits of the decimal expansion of $\pi$. Here again, mathematicians think that $\pi$ is edible.

On my website, I divided regular properties into common properties and rare properties. This separation can increase your pride if your favorite number turns out to have a rare property. There are fewer than 100 numbers below 10,000 that have one of these rare properties. At 140 below 10,000, lazy caterer and triangular numbers are not quite rare enough to be in this rare category, 
whereas pronic and square numbers, with 99 below 10,000, just make it onto the rare side.

Rare properties are: amicable, aspiring, automorphic, cake, Carmichael, Catalan, compositorial, cube, factorial, Fibonacci, Google, hungry, Mersenne, Mersenne prime, narcissistic, palindromic prime, pentagonal, perfect, power of 2, primorial, pronic, repunit, square, tetrahedral, vampire, weird.

\section{Unique Properties}

Currently I have uploaded 839 unique properties. Here are some of my favorites:
\begin{description}
\item[11] is the largest number which is not expressible as the sum of two composite numbers.
\item[19] is the largest prime which is a palindrome in Roman numerals.
\item[27] is the only number which is thrice the sum of its digits.
\item[27] is the largest number that is the sum of the digits of its cube.
\item[38] has a representation in Roman numerals --- XXXVIII --- that is lexicographically the last possible Roman numeral.
\item[40] when written in English ``forty'' is the only number whose constituent letters appear in alphabetical order.
\item[99] is the largest number that is equal to  the sum of its digits plus the product of its digits: 99 = 9 + 9 + 9 * 9.
\item[119] cents is the largest amount of money one can have in coins without being able to make change for a dollar.
\item[144] is the only composite square in the Fibonacci sequence.
\item[888] is the only number whose cube, 700227072, consists of 3 digits each occurring 3 times.
\item[1089] is the smallest number whose reverse is a non-trivial integer multiple of itself: 9801 = 9*1089.
\item[1210] is the smallest autobiographical number: $n = x_0x_1x_2\cdot x_9$ such that $x_i$ is the number of digits equal to $i$ in $n$: 1210 is built of one 0, two 1s, and one 2. Martin Gardner in his book ``Mathematical Circus'' \cite{Gardner} gives a puzzle to calculate the only 10-digit autobiographical number. See also my paper \cite{KhovanovaAuto}.
\item[1331] is the smallest non-trivial cube containing only odd digits.
\end{description}

Here are some properties that are easy to prove that you can do as an exercise:
\begin{description}
\item[5] is the only prime which is the difference of two squares of primes.
\item[6] is the only mean between a pair of twin primes which is triangular.
\item[1728] is the only compositorial cube (submitted by Sergei Bernstein).
\end{description}

I started this website for children, so my first properties were very easy and well-known like:
\begin{description}
\item[5] is the number of Platonic solids.
\item[9] is the smallest odd composite number.
\end{description}

At this point new properties have become more difficult to prove and feel more and more like research.

\section{One is the Only Triangular Cube}

I am particularly interested in numbers that simultaneously have two of my regular properties. 

For example, I wondered if there were any pronic cubes. Since pronic numbers are of the form $n(n+1)$ and cubes are $m^3$, I needed to find positive solutions to the equation 
$$n(n+1) = m^3.$$
As numbers $n$ and $n+1$ are coprime, each of them must be a cube. The only two integer cubes that differ by 1 are 0 and 1. Hence, pronic cubes do not exist.

Next I moved to triangular cubes. Obviously, 1 is both triangular and cubic. I spent a lot of my computer's processing power trying to find other triangular cubes. When I realized that 1 is the only triangular cube I could find, I started to look around for the proof. I asked my friends on the OEIS \cite{OEIS} SeqFan
mailing list if any of them knew the proof. 

I received proofs from Jaap Spies and Max Alekseyev. Here's Max's proof. Triangular numbers are of the form $n(n+1)/2$ and cubes are $m^3$. Hence, triangular cubes correspond to the solutions of the equation
$$n(n+1)/2 = m^3.$$
The big trick in solving this equation is to multiply both sides by 8:
$$4n(n+1) = 8m^3.$$
Then we can rearrange the equation into:
$$(2n+1)^2 - 1  = (2m)^3.$$
By the Catalan's conjecture, which was recently proven, the only two positive powers that differ by 1 are 9 and 8. From here the proof follows.

\section{Odd Fibonaccis}

I stumbled upon the following statement in the MathWorld \cite{MathWorld}: no odd Fibonacci is divisible by 17. So I started wondering if I can make a unique property out of it. Is 17 the only such number? Clearly not, as multiples 
of 17 will also have the same property. Is 17 the smallest such number? Obviously, no odd Fibonacci is divisible by 2. So 17 is not the smallest. What if we exclude 2? By thinking it through, I proved two new unique properties for 
my website and wrote a paper on the way \cite{KhovanovaFib}. The two properties are: 
\begin{description}
\item[9] is the smallest odd number such that no odd Fibonacci number is divisible by it.
\item[17] is the smallest odd prime such that no odd Fibonacci number is divisible by it.
\end{description}

\section{Conjectures}

Here is a sample of my unique properties that are conjectures:
\begin{description}
\item[70] is the largest known number $n$ such that $2^n$ has a digit sum of $n$.
\item[86] is conjectured to be the largest number $n$ such that $2^n$ does not contain a 0.
\item[264] is the largest known number whose square is undulating: $264^2 = 69696$.
\item[2201] is the only non-palindrome known to have a palindromic cube.
\item[3375] is the largest known cube which contains all prime digits.
\end{description}

\section{Help Me Check These Properties}

You can help me with some properties. It is not enough to send me a link to a website. I want either a proof or a program, preferably in a programming language I can read (Java, Mathematica) or a reference to a paper in a refereed journal.

I have 674 properties in my database that need checking. Here are some:
\begin{itemize}
\item The maximum number of squares a chess bishop can visit, if it is only allowed to visit each square once.
\item The number of legal knight/king/bishop moves in chess.
\item The smallest prime number that is the sum of a prime number of consecutive prime numbers in a prime number of different ways.
\item The number of forms of magic knight's tour on the chessboard.
\item The number of primes that can appear on a 24-hour and/or 12-hour digital clock.
\item The smallest number that cannot be added to a nonzero palindrome such that the sum is also palindromic.
\end{itemize}

\section{Submissions}

If you submit a property, I do not guarantee that I will add it because I already have about 1,500 properties of numbers in my database that are not true, not unique or not very interesting.

In addition, I do not like properties that contain parameters. Parameters mean that there is a sequence and therefore the property is not unique. As a rule of thumb if you can write your property without using numbers, I like it. Also, sometimes I allow 2, 3, 10 or 100 as parameters. In general, the bigger the number the more accepting and forgiving I am.

Also, if you have a proof, please send the proof, too. If I like your property and understand your proof, I will upload the property very fast and include your name as a submitter. If I do not understand the proof for your property, I 
will still have your property and your name in my internal database, but you will have to wait until I get around to proving it.

\section{Alexey Radul}

My son, Alexey Radul, suggested the idea for Number Gossip about 10 years ago and at first I uploaded it on my personal website \cite{KhovanovaWeb}. In 2006 Alexey redesigned Number Gossip using Ruby on Rails and at the end of 2006
it got its own url — http://www.numbergossip.com/. 

Alexey was also one of the first submitters. Here is one of his submissions:
\begin{description}
\item[6] is the only even evil perfect number.
\end{description}

Let me remind you that perfect numbers are numbers that are sums of their proper divisors. For example: $6 = 1 + 2 + 3$ and $28 = 1 +2 + 4 + 7 + 14$. It is not known if odd perfect numbers exist. But for even perfect numbers it is known that they are of the form: $2^{p-1}(2^p − 1)$ for prime $p$. 

Actually, we know even more than that: we know that perfect numbers are in a one-to-one correspondence with $p$, for which $2^p - 1$ is prime. That means, that even perfect numbers are in a one-to-one correspondence with Mersenne primes.

If a number is represented with the powers of two, usually it is easy to find its binary representation. From here you can prove that the evilness of 6 follows from the evenness of 2. The only even prime is 2, hence the only evil even perfect number is 6.

\section{Statistics}

Currently, the highest number you can input on  Number Gossip is 9,999. Other than the special case of number 1, every number has at least four properties, one from each of the following groups:
\begin{itemize}
\item even or odd
\item prime or composite
\item evil or odious
\item perfect, abundant or deficient
\end{itemize}

During my talk I announced that the smallest number that does not have a unique property is 32. Someone pointed out to me that this fact is in itself a unique property. I do not want to add self-referencing properties to my database, but on my way home from G4G8 I invented a property for 32:
\begin{description}
\item[32] is conjectured to be the highest power of two with all prime digits.
\end{description}

I checked this property up to $86^{th}$ power of 2. And 86 is conjectured to be the highest power of two that doesn't contain zero. That means that 32 conjecture follows from 86 conjecture.

The next number without a unique property is 51, followed by 56 and 57.
The largest number that has a unique property is currently 8833:
\begin{description}
\item[8833] is the largest 4-digit number that is the sum of the squares of its halves: $8833 = 88^2 + 33^2$.
\end{description}

\section{Future Plans}

I've got really big plans for Number Gossip:

\begin{itemize}
\item I have 647 more unique properties in my internal database that I need to check and I am always discovering new properties.
\item I have a list of about ten more regular properties I want to add, like brilliant, fortunate, primeval and totient numbers.
\item I would like to move the limit from 9,999 up to 20,000.
\item As I progress, the new properties are increasingly harder to prove. I would like to provide proofs of difficult properties on the website because I think you would find that interesting.
\item Some of my unique properties are of the form ``the smallest number that$\ldots$'' or ``the largest number that$\ldots$'' These properties can sometimes have a corresponding sequence in the online database (see OEIS \cite{OEIS}). When they do, I will be adding references to the OEIS.
\end{itemize}

\section{Acknowledgments}

I am looking for unique properties everywhere I can. I got my initial encouragement from Erich Friedman's page ``What's Special About This Number?'' (see \cite{Friedman}). Later I used the Online Encyclopedia of Integer Sequences \cite{OEIS}, Wikipedia \cite{wiki} and the ``Prime Curios!'' collection (see \cite{PrimeCurios}).

I am grateful to Sue Katz for helping me to edit this paper.

\end{document}